\documentclass[12pt]{article}
\usepackage{amsfonts}
\usepackage{xcolor}

\usepackage[plainpages,backref,urlcolor=blue]{hyperref}
\usepackage{tikz}

\title {\bf{A stable version of Terao conjecture}}

\author{\bf{Cristian Anghel}}  
\newtheorem{defn}{Definition}[section]
\newtheorem{conj}[defn]{Conjecture}
\newtheorem{tth}[defn]{Theorem}

\newtheorem{tth1}[defn]{Babylonian tower theorem}
\newtheorem{que}[defn]{Question}
\newtheorem{conj1}[defn]{Stable Terao conjecture}
\begin{document}

\maketitle

\bigskip
\noindent
{\small {{\bf ABSTRACT:} The aim of this note is to introduce a stable version of Terao conjecture, using the notion of infinitely stably extendability of 
vector bundles on $\mathbb P^n$, considered and characterized by I. Coanda in  arXiv:0907.4040. } }


\


\tableofcontents

\section{Introduction}

Terao conjecture (1981) introduced in \cite{te}, asserts that the freeness of a hyperplane arrangement depends only of its combinatorics. For a good survey see \cite{yo}. The freeness is equivalent with the fact that the associated bundle  splits completely as 
direct sum of line bundles. This last property, thanks to Horrocks criterion \cite{ho2} \cite{ho3}, is equivalent with the vanishing of certain cohomology modules of the bundle in question. Also, using the famous  Barth-Van de Ven-Sato-Tyurin result \cite{bava} \cite{sa} \cite{tyu},  the freeness of an arrangement is equivalent with the infinitely extendability of the associated bundle. In the first part of the paper, we shall describe the above circle of ideas.

The second part will be devoted to the notion of stably extendability, introduced by Horrocks (1966) in \cite{ho1}, and its connection with the above results, thanks to a theorem of Coanda (2009) from \cite{co1}, which gives a characterization of infinitely stably extendable vector bundles in terms of the vanishing of some cohomology modules of the bundle. Finally, we shall formulate a problem with the same flavor as Terao conjecture, using the Coanda notion of infinitely stably extendabiliy.

\section{Arrangements and their lattices} \label{sec:prelim}

An arrangement in the complex projective space $\mathbb P^n(\mathbb C)$ is a finite collection of 
hyperplanes ${\mathcal A}=\{ H_1,...,H_k \} .$
For a fixed arrangement ${\mathcal A}$, its intersection lattice $L_{\mathcal A}$ is the poset with elements the finite intersections 
between the $H_i's$, ordered by reverse inclusion: for $L_1,L_2\in L_{\mathcal A}$, $$L_1 \leq L_2 \ iff \ L_1 \supseteq L_2.$$
Using the interection lattice, we have a first equivalence relation for arrangements:
\begin{defn} Two arrangements ${\mathcal A}_1, {\mathcal A}_2$
have the same combinatorics if their lattices $L_{\mathcal A_1}, L_{\mathcal A_2}$ are isomorphic.
\end{defn}
For example, if ${\mathcal A}_1$ is defined by three concurrent lines in  $\mathbb P^2(\mathbb C)$ and ${\mathcal A}_2$ by three lines without a common point, then  ${\mathcal A}_1$ and ${\mathcal A}_2$ have different combinatorics.

A fundamental question  in the theory of hyperplanes arrangements is to find which properties of the arrangement depends only on its lattice i.e. only of its combinatorics.
For example, concerning the cohomology algebra of the complement we have the following celebrated result due to Arnold, Brieskorn, Orlik and Solomon \cite{ar} 
\cite{br} \cite{os}:
\begin{tth}
The cohomology ring $H^*(\mathbb P^n(\mathbb C) \ \backslash \ \bigcup\limits_{i=1}^{k} H_i)$ of the complement of ${\mathcal A}=\{ H_1,...,H_k \} $ is combinatorially determined by $L_{\mathcal A}$.
\end{tth}
Also, a negative result in this direction, concern the homotopy type of the complement: for example $\pi _1(\mathbb P^n(\mathbb C) \ \backslash \ \bigcup\limits_{i=1}^{k} H_i)$ is not combinatorially determined.
In fact, Rybnikov (1998) in \cite{ry} constructed two arrangements in $\mathbb P^2(\mathbb C)$ with the same combinatorics but different $\pi _1$ for the complements.

\section{Bundles associated with arrangements } \label{sec:bogo}

Apart the lattice and homological or homotopical invariants associated to an arrangement ${\mathcal A}$ another interesting object is the sheaf ${\mathcal T}_{\mathcal A}$ of vector fields with logarithmic poles along ${\mathcal A}$. It was introduced for the first time by Saito and Deligne in the '80s \cite{sai} and used in the context of hyperplane arrangements by Dolgachev, Kapranov \cite{do1} \cite{do2}, Terao and others. Its construction goes as follows:
denote by $f_i$ an homogenous equation of the hyperplane $H_i$ and by $f$ the product $\prod_{i=1}^{k}f_i$.
Then ${\mathcal T}_{\mathcal A}$ is defined as the kernel of the  map 
$${\mathcal O}_{\mathbb P^n}^{\oplus (n+1)} \rightarrow {\mathcal O}_{\mathbb P^n}(k-1),$$
defined by the partial derivatives of $f$: $(\partial_{x_0}f,...,\partial_{x_n}f)$.
The sheaf ${\mathcal T}_{\mathcal A}$ will be the principal object of study in the sequel. In general it is a rank-$n$ sheaf on $\mathbb P^n$, but we will be interested mainly in the case where it is locally free. For example, due to a result of Dolgachev this is the case for the normal crossing arrangements.

An important problem concerning ${\mathcal T}_{\mathcal A}$, in the case when it is locally free, is its splittability:
\begin{defn}
A vector bundle on $\mathbb P^n$ is splittable if it is direct sum of line bundles.
\end{defn}
With the definition above, an arrangement ${\mathcal A}$ is called free if the associated sheaf ${\mathcal T}_{\mathcal A}$ is splittable.
Of course if ${\mathcal A}$ is free, then ${\mathcal T}_{\mathcal A}$ is locally free and consequently, concerning the freeness one can consider only arrangements with locally free ${\mathcal T}_{\mathcal A}$. 
In the above terms, one can enounce the Terao's conjecture:
\begin{conj}
The freeness of an arrangement is combinatorially determined. Namely, for two arrangements ${\mathcal A}_1,{\mathcal A}_2$ with isomorphic lattices, if 
 ${\mathcal A}_1$ is free then ${\mathcal A}_2$ is also free.
\end{conj}

\section{Freeness versus infinitely extendability } \label{sec:free}

As long as the freeness of ${\mathcal A}$ means the splittability of ${\mathcal T}_{\mathcal A}$, a good starting point in the study of free 
arrangements could be a criterion which ensure the splittability of a vector bundle on $\mathbb P^n$.
In this direction, the fundamental result is Horrocks theorem. Let $F$ a vector bundle on $\mathbb P^n$. For any $1\leq i\leq n-1$ we denote by $H^{i}_*(F)$ the cohomology module $$\bigoplus\limits_{k\in \mathbb Z}H^i(\mathbb P^n, F(k)),$$
where, as usual $F(k)=F\otimes {\mathcal O}_{\mathbb P^n}(k)$.
With the above notations we have the following criterion of Horrocks:
\begin{tth}
A vector bundle $F$ on $\mathbb P^n$ splitts completely as direct sum of line bundles iff for any $1\leq i\leq n-1$ the cohomology module $H^{i}_*(F)$ is zero.
\end{tth}
Consequently, the freeness of an arrangement ${\mathcal A}$ is equivalent with the vanishing of all intermediate cohomology modules $H^{i}_*({\mathcal T}_{\mathcal A})$ of its bundle of logarithmic vector fields.

\

\noindent A second viewpoint concerning the splittability of bundles on $\mathbb P^n$ is connected with the following definition: 
\begin{defn}
A vector bundle $F$ on $\mathbb P^n$ is infinitely extendable if for any $m\geq n$ there exist a bundle $F_m$ on $\mathbb P^m$ such that 
$${F_m}_{\mid \mathbb P^n} \ \simeq  \ F.$$
\end{defn}
The following result, due to Barth, Van de Ven, Sato and Tyurin, asserts that in fact the infinitely extendability is equivalent with the complete splittability of the bundle in question:
\begin{tth1}
For a vector bundle $F$ on $\mathbb P^n$, the following are equivalent:\\
1. $F$ splitts completely as direct sum of line bundles,\\
2. $F$ is infinitely extendable.
\end{tth1}
As consequence, one obtain another characterization of the freeness of an arrangement ${\mathcal A}$, namely the infinitely extendability of ${\mathcal T}_{\mathcal A}$.

\

\noindent The conclusion of the above results is that the freeness of an arrangement, which is the main property in the statement of the Terao conjecture, admits at least two equivalent formulations:

\

\begin{minipage}{4cm} \textcolor{yellow}{vanishing of all \\ the intermediate \\ cohomology  of ${\mathcal T}_{\mathcal A}$} \end{minipage}
\begin{minipage}{4cm} $\Leftrightarrow $ \ \ \ \ \textcolor{red}{freeness} \ \ $\Leftrightarrow $  \end{minipage}
\begin{minipage}{4cm} \textcolor{blue}{infinitely \\extendability of ${\mathcal T}_{\mathcal A}$}. \end{minipage}

\

\noindent The main question we will discus in the sequel is the following:
\begin{que}
Is there a weaker (than freeness) property with a similar cohomological and geometrical flavor which could be used in a modified form of Terao conjecture?
\end{que}
The answer is yes and is connected with the notion of infinitely stably extendability, characterized by Coanda in 2009.

\section{The stably freeness of arrangements } \label{sec:stabfree}

In analogy with the previous notion of extendability, Horrocks (1966) introduced the following weaker concept:
\begin{defn}
A vector bundle $F$ on $\mathbb P^n$ is stably extendable on a larger space $\mathbb P^m$ if there exists a bundle $F_m$ on $\mathbb P^m$ whose restriction to 
$\mathbb P^n$ is the direct sum between $F$ and certain line bundles.
\end{defn}

\noindent A first remark is that an extendable bundle is obviously stably extendable, but the converse is not true. For example the tangent bundle of $\mathbb P^n$, $T_{\mathbb P^n}$ is stably extendable but not extendable.
Also, one should note that the above notion is connected with the complete splittability of a bundle by the following result of Horrocks:
\begin{tth}
If the bundle $F$ on $\mathbb P^n$ extends stably to $\mathbb P^{2n-3}$ and the cohomology modules $H^{1}_*(F)$, $H^{n-1}_*(F)$ vanishes, then $F$ splits completely on $\mathbb P^n$.
\end{tth} 

\noindent The result above, show that the condition of stably extendability of a bundle $F$, has a subtle connection with the property of complete splittability and is also a good motivation for the following definition introduced (and as we shall see, characterized) by Coanda in 2009:
\begin{defn}
A vector bundle $F$ on $\mathbb P^n$ is infinitely stably extendable, if for any $m\geq n$ it extends stably on $\mathbb P^m$.
\end{defn}

\noindent As in the case of stably extendability, the above property is strictly weaker than infinitely extendability, as long as again, the example of the tangent bundle of $\mathbb P^n$ shows that there are bundles infinitely stably extendable which are not splittable and therefore (using the babylonian tower theorem of Barth-Van de Ven-Sato-Tyurin) are not infinitely extendable.

The main point concerning the above property is that, like the complete splittability and therefore -via the babylonian tower theorem- like the infinitely extendability, it admits an analogous cohomological characterization in terms of some intermediate cohomology modules. This one, was obtained by Coanda in 2009:
\begin{tth}
A vector bundle $F$ on $\mathbb P^n$ is infinitely stably extendable iff for any $2\leq i\leq n-2$ the intermediate cohomology module $H^{i}_*(F)$ vanishes.
\end{tth}

\noindent On should remark that in fact, the original theorem of Coanda, contains also a third characterization of the infinitely stably extendability, namely as the condition for $F$ of being the cohomology of a free monad. However we do not use this fact for the moment. 
Also, one should note that the condition in the theorem is empty for $n\leq 3$ but it can be proved that any bundle on $\mathbb P^{\leq 3}$ is infinitely stably extendable \cite{co2}.

\section{A "stable" Terao conjecture } \label{sec: stabterao}

Inspired by the above result we consider de following:
\begin{defn}
An arrangement ${\mathcal A}$ is stably free if its associated bundle of vector fields with logarithmic poles ${\mathcal T}_{\mathcal A}$ is infinitely stably extendable.
\end{defn}

\noindent Obviously, a free arrangement is stably free, the converse is not true and we have, as in the case of freeness,  a similar characterization of stably freeness:

\

\begin{minipage}{4cm} \textcolor{yellow}{vanishing of all \\ the intermediate \\ cohomology  of ${\mathcal T}_{\mathcal A}$ in the range 
$2\leq i\leq n-2$}\end{minipage}
\begin{minipage}{4.5cm} $\Leftrightarrow $ \ \ \textcolor{red}{stably freeness} \ \ $\Leftrightarrow $ \ \ \  \end{minipage}
\begin{minipage}{4cm} \textcolor{blue}{infinitely stably \\extendability of ${\mathcal T}_{\mathcal A}$}. \end{minipage}

\

\

\noindent Consequently, we introduce the following:
\begin{conj1}
The stably freeness of an arrangement is \\ combinatorially determined. Namely, for two arrangements ${\mathcal A}_1,{\mathcal A}_2$ with isomorphic lattices, if 
 ${\mathcal A}_1$ is stably free then ${\mathcal A}_2$ is also stably free.
\end{conj1}

\noindent A first remark is that the two conjectures are not comparable: no one implies the other.
Also, it could be interesting to compare this notion of stably freeness which can be obviously be extended from arrangements to arbitrary union of hyper-surfaces to other weaker notions of freeness existing in literature.
For convenience we mention only two:\\
- the nearly free and almost free divisors introduced by Dimca and Sticlaru \cite{dist},\\
- the quasi free divisors studied by Castro-Jimenez \cite{uch}. 

\

\noindent Another point to note is the following. Due to Horrocks and Coanda criteria, both conjectures can be expressed as 
the  combinatorial invariance of the vanishing in a certain range of the intermediate cohomology modules of ${\mathcal T}_{\mathcal A}$. From this viewpoint, one can ask the following question which already appeared in \cite{yo} in relation with the lattice cohomolgy introduced by Yuzvinsky in \cite{yu1} \cite{yu2} \cite{yu3}:

\begin{que}
For an arrangement ${\mathcal A}$ in $\mathbb P^n$, and a fixed $1\leq i\leq n-1$, can be the cohomology module $H^{i}_*({\mathcal T}_{\mathcal A})$ expressed/computed only in terms of the lattice  $L_{\mathcal A}$ of ${\mathcal A}$ ?
\end{que}

\bigskip

\noindent
Cristian Anghel \\
Department of Mathematics\\
Institute of Mathematics of the Romanian Academy\\
Calea Grivitei nr. 21 Bucuresti Romania\\
email:\textit{Cristian.Anghel@imar.ro}

\end{document}